\documentclass[12pt,amsmath,amscd,amsbsy,amssymb]{amsart}

\theoremstyle{plain}
\newtheorem{thm}{Theorem}[section]
\newtheorem{lemma}[thm]{Lemma}

\newtheorem{prop}[thm]{Proposition}

\theoremstyle{definition}
\newtheorem{defin}[thm]{Definition}

\newcommand{\N}{{\mathbb N}}
\newcommand{\F}{{\mathcal F}}
\newcommand{\G}{\mathcal G}
\newcommand{\Oo}{\mathcal O}
\newcommand{\VG}{{\mathcal V}{\mathcal G}}
\newcommand{\EG}{{\mathcal E}{\mathcal G}}
\newcommand{\comp}{{\mathrm C}{\mathrm o}{\mathrm m}{\mathrm p} }
\newcommand{\Cc}{{\mathcal C}}
\newcommand{\Ss}{{\mathcal S}}

\newcommand{\pr}{{\mathrm p}{\mathrm r}}
\numberwithin{equation}{section}

\usepackage[all]{xy}

\begin{document}

\title[]
{Open-multicommutativity of the functor of upper-continuous
capacities}
\author{Roman Kozhan}

\address{Department of Mechanics and Mathematics, Lviv National University,
Universytetska 1, 79000 Lviv, Ukraine}

\email{}

\thanks{The author gratefully thanks Michael Zarichnyi for helpful ideas,
 discussions and comments.}

\begin{abstract}The notion of open-multicommutativity, introduced by
Kozhan and Zarichnyi \cite{KozZar}, is investigated. The weakly
normal covariant functor of upper-continuous capacities is
considered. The main result of the paper is that this functor
open-multicommutative.
\end{abstract}

\maketitle


\section{Introduction}\label{s:intro}

The impact of the non-additive probability theory on the modern
economics and finance increased significantly during the last
decades. This theory is based on the notion of capacity (also known
as non-additive measure) which was first introduced by Choquet
\cite{Choq}. By the 80's the number of authors (Schmeidler
\cite{Schmeidl}, Quiggin \cite{Quig}, Yaari \cite{Yaari}) presented
axiomatizations of individual's preferences and developed the
non-expected utility theory which based on the notion of the Choquet
integral.

From the topological point of view capacities were considered by
Zhou \cite{Zhou}. He investigated the structure of the space of
upper-continuous capacities and established an integral
representation of continuous comonotonically additive functionals.

In this paper we study the space of upper-continuous capacities from
the viewpoint of the categorical topology. We prove an analogical
result which was investigated in the case of the probability
measures space. The notion of open-multicommutativity which combines
properties of a covariant functor to be open and bicommutative has
been introduced in \cite{KozZar}. The main result of their paper is
that the functor of probability measures is open-multicommutative in
the category of compact Hausdorff spaces. Here we extend an area of
this investigation and consider the functor of upper-continuous
capacities. Although, this functor turns out to be weakly-normal, it
satisfies the open-multicommutativity property.

The paper is organized as follows. In Section \ref{s:part2} we
remind some definitions which we use below. Section \ref{s:part3}
contains a proof of the finite open-multicommutativity of the
capacity functor. The main result is given at the end of the last
section.

\section{Notations and Definitions}\label{s:part2}

\subsection{Functor of upper-continuous capacities.}

Let $X$ be a compact Hausdorff space and $\F$ a $\sigma$-algebra
of its Borel subsets.

\begin{defin}\label{d:capac} A real-valued set function $\mu$ on $\F$ is
called a {\em capacity} if $\mu(\emptyset)=0$, $\mu(X)=1$ and
$\mu(A)\leq\mu(B)$ for all $A\subseteq B$, $A,B\in\F$.
\end{defin}

\begin{defin}\label{d:upper} A capacity $\mu$ is {\em upper-continuous} if
$\underset{n\rightarrow
\infty}{\lim}\mu(A_n)=\mu(\underset{n=1}{\overset{\infty}{\cap}}A_n)$
for any monotonic sequence of sets $A_1\supseteq A_2\supseteq
A_3\supseteq...$ with $A_n\in\F$, $n\in\N$.
\end{defin}

We denote a set of all upper-continuous capacities on $X$ as $M(X)$.
Due to Zhou \cite{Zhou} we can identify the set $M(X)$ with the set
of all comonotonically additive, monotonic and continuous functional
on $C(X)$ by the formula
\[\mu(f)=\underset{0}{\overset{\infty}{\int}}\mu(f\geq t)dt+
\underset{-\infty}{\overset{0}{\int}}(\mu(f\geq t)-1)dt.
\]The above integral is called the Choquet integral.

Let us endow the set $M(X)$ with the weak-* topology. The base of
this topology consists of the set of the form
\[O(\mu_0, f_1,...,f_n, \varepsilon)=
\{\mu\in M(X)\colon |\mu_0(f_i)-\mu(f_i)|<\varepsilon,
i=1,...,n\},
\]where $\mu_0\in M(X)$, $f_1,...,f_n\in C(X)$ and $\varepsilon>0$.

We can consider the map $M\colon\comp\to\comp$ as a covariant
functor in the category $\comp$ and as it is shown in \cite{NykZar}
that it is also weakly normal. Another important property of this
functor is that it is open and bicommutative.

\begin{prop}\label{l:1} Functor $M$ is bicommutative.
\end{prop}

\begin{proof}Let us consider an arbitrary bicommutative diagram
\begin{equation}\label{e:sqdiag}\xymatrix{Z\ar[r]^f\ar[d]_g &X\ar[d]^h\\Y\ar[r]_s &T}
\end{equation}in the category $\comp$. In order to prove that $M$ is
a bicommutative functor it is sufficient to show that for every
$\mu\in M(X)$ and $\nu\in M(Y)$ such that $Mh(\mu)=Ms(\nu)=\tau\in
M(T)$ there exists a capacity $\lambda\in M(Z)$ with

\begin{equation}\label{e:marg}Mf(\lambda)=\mu \text{ and } Mg(\lambda)=\nu.
\end{equation}
Due to condition (\ref{e:marg}) for every $A\in\F_X$ and
$B\in\F_Y$ it must hold
\[\lambda(f^{-1}(A))=\mu(A)\text{ and }\lambda(g^{-1}(B))=\nu(B).
\]Denote $\Ss=\{f^{-1}(A),g^{-1}(B)\colon A\in\F_X,B\in\F_Y\}$. Let
$\lambda$ be an inner measure defined as
\[\lambda(D)=\sup\{\mu(f(C)),\nu(g(C))\colon C\subseteq D, C\in\Ss\}
\]for every $D\in\F_Z$. Defined in such way set function $\lambda$ is an upper-continuous
capacity (see, for instance, \cite{Zhang}). Let us show that
condition (\ref{e:marg}) is satisfied. Let $A\in\F_X$ and
$A^\prime=f^{-1}(A)\subset Z$. Obviously that $\sup\{\mu(f(C))\colon
C\subseteq A^\prime, C\in\Ss\}=\mu(A)$. We assume that there exists
a subset $B\subset Y$ such that $B^\prime=g^{-1}(B)\subseteq
A^\prime$ and $\nu(B)>\mu(A)$. Note that the set
$\tilde{A}=h^{-1}(s(B))\subseteq A$. Indeed, due to the definition
of this set for every point $a\in\tilde{A}$ we can find $b\in B$
such that $h(a)=s(b)$. Because of the bicommutativity of diagram
(\ref{e:sqdiag}) there exists point $z\in Z$ satisfying $f(z)=a$ and
$g(z)=b$. Since $B^\prime$ is a full preimage of the set $B$ it is
necessary that $b\in B^\prime\subseteq A^\prime$. This implies $a\in
A$. Due to the condition (\ref{e:marg}) we have
\[\nu(B)\leq
\tau(s(B))=\mu(h^{-1}(s(B)))=\mu(\tilde{A})\leq\mu(A),
\]which contradicts our assumption. Thus,
$\lambda(f^{-1}(A))=\mu(A)$ for every $A\in \F_X$. Analogically we
can prove that $\lambda(g^{-1}(B))=\nu(B)$ for each $B\in\F_Y$.
Therefore, condition (\ref{e:marg}) is satisfied.
\end{proof}

\begin{prop}\label{l:2} Functor $M$ is open.
\end{prop}

\begin{proof} This proposition is proved in \cite{NykZar}.
\end{proof}

\subsection{Open-multicommutative functors and characteristic map.}

Let us recall the notion of the multi-commutativity of a
weakly-normal functor which is introduced in \cite{KozZar}.

Suppose that $\G$ is a finite partially ordered set and we also
regard it as a finite directed graph. Denote by $\VG$ the class of
all vertices of graph $\G$ and by $\EG$ the set of its edges. A
functor $\Oo\colon\G\to\comp$ is called a {\em diagram}. A {\em
cone} over $\Oo$ consists of a space $X\in|\comp|$ and a family of
maps $\{X\to\Oo(o)\}_{o\in\VG}$ that satisfy obvious commutativity
conditions. Given such a cone, $\Cc=(\{X\to\Oo(o)\}_{o\in\VG})$, we
denote by $\chi_{\Cc}\colon X\to\lim\Oo$ its characteristic map.

We say that the cone $\Cc$ is {\it open-multicommutative} if its
characteristic map is an open onto map.

\begin{defin}A normal functor $F$ in $\comp$ is called {\em
open-multicommutative} ({\em finite open-multicommutative}) if it
preserves the class of open-multicommutative diagrams (which consist
of finite spaces).
\end{defin}

The following result can be found in \cite{Koz}.

\begin{prop}\label{p:1} For a weakly normal open bicommutative functor $F$
the following properties are equivalent:

$F$ is open-multicommutative;

$F$ is finite open-multicommutative.
\end{prop}

\section{Open-multicommutativity of $M$}\label{s:part3}

Let us assume first that all spaces $\Oo(o)$, $o\in\VG$ are finite
and discrete. According to Proposition \ref{p:1} for the
open-multicommutativity of the functor $M$ it is sufficient to show
that it is finite open-multicommutative, i.e. the characteristic map
$\chi\colon M(\lim\Oo)\to\lim M(\Oo)$ is open and surjective.

Let us also recall that $\lim\Oo$ can be defined in terms of {\em
threads}. We say that the point $x=(x_o)_{o\in\VG}\in
\underset{o\in\VG}\prod\Oo(o)$ is a {\it thread} of the diagram
$\Oo$ if for every $o_1,o_2\in\VG$ with $o_1\leq o_2$ it holds
$\pr_{o_1}(x)=\varphi_{o_1o_2}\circ\pr_{o_2}(x)$. It is well known
that $\lim\Oo\subseteq\underset{o\in\VG}{\prod}\Oo(o)$ and since all
$\Oo(o)$ are discrete, the limit of the diagram is also discrete
space.

Let $\lambda^0\in M(\lim\Oo)$ be a capacity on the space $\lim\Oo$
and $\mu^0_o\in M(\Oo(o))$ be its marginals for $o\in\VG$. Let
$O(\lambda^0, f_1,...,f_n,\varepsilon)$ be an arbitrary weak-*
neighborhood of the non-additive measure $\lambda^0$.

In order to prove the openness of the characteristic map it is
sufficient to find a neighborhood of $(\mu^0_o)_{o\in\VG}$ such that
every point from this neighborhood can be covered by some capacity
from $O(\lambda^0, f_1,...,f_n,\varepsilon)$.

\begin{lemma}\label{l:3} Let $X$ be a discrete compactum. The base
of the weak-* topology on $M(X)$ consists of the sets of the form
\[O(\mu, F_1,..., F_n, \varepsilon)=
\{\nu\in M(X)\colon|\nu(F_i)-\mu(F_i)|<\varepsilon, i=1,...,n \},
\]$\mu\in M(X)$ and $F_i\subset X$, $i=1,...,n$.
\end{lemma}

\begin{proof} Let us show first that for every set of the form
$O(\mu, F_1,...,F_n, \varepsilon)$ we can find a basis
neighborhood $O(\mu, f_1,...,f_k, \delta)$ for some $f_i\in C(X)$
and $\delta>0$, $i=1,...,k$. Indeed, if we set $f_i={\bf 1}_{F_i}$
and $\delta=\varepsilon$ we get
\[O(\mu, F_1,..., F_n)=O(\mu, f_1,...,f_n,\varepsilon).
\]Conversely, consider an element of the sub-base $O(\mu,
f,\varepsilon)$. Since the space $X$ is discrete we can represent
$f=\underset{i=1}{\overset{k}\sum}\alpha_i {\bf 1}_{F_i}$ such that
$F_1\subset F_2\subset...\subset F_k$. It is clear (see \cite{Ang})
that for every capacity $\nu\in M(X)$ it holds
$\nu(f)=\underset{i=1}{\overset{k}\sum}\alpha_i\nu(F_i)$. Let us
consider a set $O(\mu, F_1,.., F_k, \frac{\varepsilon}{k\alpha})$,
where $\alpha=\max\{|\alpha_1|,...,|\alpha_k|\}$. Comonotonicity of
functions ${\bf 1}_{F_i}$ implies that for every capacity $\nu\in
O(\mu, F_1,.., F_k, \frac{\varepsilon}{k\alpha})$ we have
\[|\mu(f)-\nu(f)|=|\underset{i=1}{\overset{k}\sum}\alpha_i
(\mu(F_i)-\nu(F_i))|\leq
\underset{i=1}{\overset{k}\sum}|\alpha_i||\mu(F_i)-\nu(F_i)|
<\underset{i=1}{\overset{k}\sum}\frac{\varepsilon\alpha_i}{k\alpha}<\varepsilon.
\]Hence, $O(\mu, F_1,..,
F_k, \frac{\varepsilon}{k\alpha})\subset O(\mu, f,\varepsilon)$.
\end{proof}

According to Lemma \ref{l:3} we can assume without loss of
generality that functions $f_1,...,f_n$ are of the form
$f_i(x)=\left\{\begin{array}{ll}1,&x\in F_i,\\0,&x\notin
F_i\end{array}\right.$ for every $x\in \lim D$ and some
$F_1,...,F_n\subseteq\lim\Oo$.

We consider a neighborhood
\[U=O(\mu^0_1,\{x^1_1\},...,\{x^1_{m_1}\},\delta)\times...\times
O(\mu^0_k,\{x^k_1\},...,\{x^k_{m_k}\},\delta),
\]where
$X_i=\{x^i_1,...,x^i_{m_i}\}$. Let $(\mu_1,...,\mu_k)$ be arbitrary
point in $U$. Let us define a capacity $\lambda$ on $\lim\Oo$.

For every subset $A\subset\lim\Oo$ we denote
\[l_A=\underset{o^\prime\in\VG}{\max}\{\max\{\mu_{o^\prime}(W)
\colon W\subseteq\Oo(o^\prime),
(\underset{o\in\VG\setminus\{o^\prime\}}{\prod}\Oo(o)\times
W)\cap\lim\Oo\subseteq A\}\}.
\]Analogically,
\[u_A=\underset{o^\prime\in\VG}{\min}\{\min\{\mu_{o^\prime}(W)
\colon W\subseteq\Oo(o^\prime),
(\underset{o\in\VG\setminus\{o^\prime\}}{\prod}\Oo(o)\times
W)\cap\lim\Oo\supseteq A\}\}.
\]Note that the interval $[l_A,u_A]$ is not empty and in order
$\lambda$ to be well defined it should satisfies inequalities
\[l_A\leq \lambda(A)\leq u_A
\]for every subset $A\subset\lim\Oo$. Recall also that
$(\mu_1,...,\mu_k)\in U$ and this implies that for every
$A\subset\lim\Oo$ we have $l_A-\delta<\lambda^0(A)<u_A+\delta$.

\begin{lemma}\label{l:4} If $A\subseteq B$ then $l_A\leq l_B$ and
$u_A\leq u_B$.
\end{lemma}

\begin{proof} It is clear that for every $W\subseteq\Oo(o^\prime)$ such
that $(\underset{o\in\VG\setminus\{o^\prime\}}{\prod}\Oo(o)\times
W)\cap\lim\Oo\subseteq A$ we have that
$(\underset{o\in\VG\setminus\{o^\prime\}}{\prod}\Oo(o)\times
W)\cap\lim\Oo\subseteq B$ for every $j=1,...,k$. This implies that
$l_A\leq l_B$.

The analogical result for the upper bounds can be derived from the
statement
\[(\underset{o\in\VG\setminus\{o^\prime\}}
{\prod}\Oo(o)\times W)\cap\lim\Oo\supseteq A\supseteq B.
\]\end{proof}

For a Borel set $A\subset \lim\Oo$ we set
\[\lambda(A)=\max\{l_A,\min\{u_A,\lambda^0(A)\}\}.
\]

\begin{lemma}\label{l:5} The set function $\lambda$ is a well-defined
capacity.
\end{lemma}

\begin{proof} First of all, $l_{(\lim\Oo)}= u_{(\lim\Oo)}=1$ this
implies that $\lambda(\lim\Oo)=1$.

$l_{\emptyset}=u_{\emptyset}=0$ this implies that
$\lambda(\emptyset)=0$.

Let us check now a monotonicity of $\lambda$. We suppose that
$A\subset B\subset \lim\Oo$. Consider the following three cases:

1). $\lambda^0(A)\in[l_A,u_A]$. In this case
\[l_A<\lambda(A)=\lambda^0(A)\leq\min\{u_A,\lambda^0(B)\}
\leq\min\{u_B,\lambda^0(B)\}=\lambda(B).
\]

2). $\lambda^0(A)>u_A$. We have
\[l_A<\lambda(A)=u_A\leq\min\{u_B,\lambda^0(A)\}\leq
\min\{u_B,\lambda^0(B)\}=\lambda(B).
\]

3). $\lambda^0(A)<l_A$. This condition implies that
\[\lambda(A)=l_A\leq l_B\leq\lambda(B).
\]\end{proof}

Let us set now $\delta=\varepsilon$. In this case we obtain for
every $i=1,...,n$ the relationship
\[|\lambda(F_i)-\lambda^0(F_i)|<\delta=\varepsilon.
\]This leads to $\lambda\in O(\lambda^0,f_1,...,f_n,\varepsilon)$.

Due to the definition of $l_A$ and $u_A$ it is easy to check that
\[l_{((\underset{o\in\VG\setminus\{o^\prime\}}{\prod}\Oo(o)\times
W)\cap\lim\Oo)}=\mu_{o^\prime}(W)
\]and
\[u_{((\underset{o\in\VG\setminus\{o^\prime\}} {\prod}\Oo(o)\times
W)\cap\lim\Oo)}=\mu_{o^\prime}(W)
\]for every $o^\prime\in\VG$ and $W\subset X_{o^\prime}$. This implies that
$\lambda((\underset{o\in\VG\setminus\{o^\prime\}}
{\prod}\Oo(o)\times W)\cap\lim\Oo)=\mu_{o^\prime}(W)$ and hence
$M\pr_{o^\prime}(\lambda)=\mu_{o^\prime}$ for all
$o^\prime\in\VG$.

Hence we proved that the inverse to the correspondence map is open
in the case of discrete $\Oo(o)$, $o\in\VG$. Thus, applying this
fact to Proposition \ref{p:1} we obtain

\begin{thm}\label{t:2} The correspondence map $\chi$ of the
diagram $\Oo$ is open and surjective for every $\Oo(o)\in|\comp|$,
$o\in\VG$.
\end{thm}

A special case of open-multicommutativity was considered by Eifler
\cite{Eifler}. One can get this case setting the set
$\EG=\emptyset$. L. Eifler proved that the functor of the
probability measures preserves surjectivity and openness of the
characteristic maps of such kind of diagrams. Thus, the result of
Theorem \ref{t:2} is an extension of Eifler's theorem on the case of
non-additive measures.

\end{document}